\documentstyle[12pt,amsfonts,amscd]{amsart}
\newtheorem{Theorem}{Theorem}[section]
\newtheorem{Definition}[Theorem]{Definition}
\newtheorem{Proposition}[Theorem]{Proposition}

\newtheorem{Lemma}[Theorem]{Lemma}
\newtheorem{Corollary}[Theorem]{Corollary}
\theoremstyle{remark}
\newtheorem{Example}[Theorem]{Example}

\def\dim{\operatorname{dim}}

\def\be{\begin{enumerate}}
\def\ee{\end{enumerate}}
\def\bT{\begin{Theorem}}
\def\eT{\end{Theorem}}
\def\bP{\begin{Proposition}}
\def\eP{\end{Proposition}}
\def\bD{\begin{Definition}}
\def\eD{\end{Definition}}
\def\bE{\begin{Example}}
\def\eE{\end{Example}}
\def\bL{\begin{Lemma}}
\def\eL{\end{Lemma}}
\def\bC{\begin{Corollary}}
\def\eC{\end{Corollary}}

\begin{document}
\title{Isolated fixed point sets for holomorphic maps}
\author{Buma L. Fridman, Daowei Ma and Jean-Pierre Vigu\'{e}}
\begin{abstract}
We study discrete fixed point sets of holomorphic self-maps of complex manifolds. The main attention is focused on the cardinality of this set
and its configuration. As a consequence of one of our observations, a
bounded domain in ${\Bbb C}^n$ with no non-trivial holomorphic retractions is
constructed.

\smallskip
\font\sc=cmcsc10
\noindent
{\sc R\'esum\'e}. Nous \'etudions les ensembles discrets qui sont les ensembles
 de points fixes d'une application holomorphe d'une vari\'et\'e complexe
 dans elle-m\^eme. En particulier, nous \'etudions le nombre d'\'el\'ements
 de ces ensembles et leurs configurations. Comme application de ces
 r\'esultats, nous contruisons un domaine born\'e de ${\Bbb C}^n$ sans r\'etraction
 holomorphe non triviale. 

\end{abstract}
\keywords{}
\subjclass[2000]{Primary: 32M05, 58C30}
\address{ fridman@@math.wichita.edu, Department of Mathematics,
Wichita State University, Wichita, KS 67260-0033, USA}
\address{ dma@@math.wichita.edu, Department of Mathematics,
Wichita State University, Wichita, KS 67260-0033, USA}
\address{ vigue@@math.univ-poitiers.fr, UMR CNRS 6086, Universit\'{e} de Poitiers, 
Math\'{e}matiques, SP2MI, BP 30179, 86962 FUTUROSCOPE, FRANCE}
\maketitle \setcounter{section}{-1}
\section{Introduction}

In classical mechanics the following Euler's theorem is well known: the
general displacement of a rigid body with one point fixed is a rotation
about some axis. So, if one considers an orientation-preserving isometry of a domain in ${\Bbb R}^3$
fixing one point, the fixed point set of this isometry will necessarily
contain at least a segment, so the fixed point set cannot be a discrete set. In the euclidean space ${\Bbb R}^n$,
one can always find a domain which has a euclidean isometry with exactly {\it one}
fixed point, however for any $n$, if an isometry of a domain in ${\Bbb R}^n$
has {\it two} fixed points it will force the existence of at least a segment
to belong to the fixed point set, and so this set will be at least one dimensional.

Switching to complex analysis, we remark that any holomorphic automorphism
of a bounded domain in ${\Bbb C}^n$ (or in general, hyperbolic manifold) is
an isometry in an invariant metric, so an Euler type statement is certainly
meaningful, that is if this automorphism has a discrete fixed point set one
can inquire what its cardinality and structure might be. To
describe this more precisely, let $f:M\rightarrow M$ be a holomorphic
self-map of a complex manifold $M$. Let $Fix(f)$ denote the set of fixed
points $\{x\in M$ $|$ $f(x)=x\}$ of $f$. Suppose that this set is discrete.
In this paper we shall examine mostly two questions. First, how large this
set can be for specific cases: $M$ is a bounded domain in ${\Bbb C}^n$, a
hyperbolic manifold, etc., while $f$ is a holomorphic automorphism or
endomorphism. Second, the structure of $Fix(f)$, namely which points of $M$
could form such a set for some holomorphic self-map of $M$. Everywhere below
we consider only holomorphic self-maps (automorphisms or endomorphisms) of
various complex manifolds, and for the sake of compactness the word {\it 
holomorphic} may be omitted. 

In examining the cardinality of a discrete fixed point set, let's first
consider the situation in one dimension. For a bounded domain $D\subset
\subset {\Bbb C}$ the discrete fixed point set of a holomorphic map $
f:D\rightarrow D$ can have no more than two points. This follows from the
following observation: any map fixing at least two points must be an
automorphism (H. Cartan), and any automorphism fixing three points must be
the identity [PL]. An annulus gives an example of a domain that has an
automorphism with exactly two fixed points.

In ${\Bbb C}^n$ the situation is not yet completely clear. For a convex
domain one has an Euler type theorem: the isolated fixed point set of any
endomorphism consists of at most one point (see Prop. 1.1). For a bounded
strictly pseudoconvex domain in ${\Bbb C}^n$ with real analytic boundary the number of points in a
discrete fixed point set of an automorphism is finite (see Thm \ref{bd}).
Must the cardinality of an isolated fixed point set of an automorphism or
endomorphism be bounded by a number depending only on the dimension of the
manifold under consideration? As one can see below (in section 1) for
endomorphisms of bounded domains in ${\Bbb C}^n$ the answer is negative. It
is also negative for automorphisms of a general hyperbolic manifold and the
entire ${\Bbb C}^n$. However, for an automorphism of a bounded domain in $%
{\Bbb C}^n$ the answer is not yet clear (see more discussions on that in the
collection of problems at the conclusion of this paper: section 5). 

We then turn to the consideration of which {\it single} points of a domain
can form $Fix(f)$ for a holomorphic $f$. 
In section 2 we prove that if every point of a hyperbolic manifold is $Fix(f)$ for some automorphism $f$ then the manifold must
be homogeneous; we also show an example of a one dimensional non-homogeneous domain with infinite number of such points.
In section 3 we consider {\it pairs of points} as fixed point sets, and
prove that for any domain ``most'' pairs of points, if fixed, force a whole
analytic set of complex dimension one to be fixed (compare to an Euler type
statement above for a domain in ${\Bbb R}^n$). 

An application is given in section 4: we construct a bounded domain in ${\Bbb C}^n$ such that if any holomorphic endomorphism of it fixes two distinct points, it will necessarily be the identity. As a consequence, this domain will have no non-trivial holomorphic retractions.

\section{General cardinality statements.}

Below we use the following notation: if $f: M\rightarrow M$ is a holomorphic self-map of $M$, then $Fix(f)$ is its fixed point set, and if such a set is discrete then $\#(Fix(f))$ is its cardinality.
We start with two positive statements (e.g. the cardinality $\#(Fix(f))$ is bounded).
\bP Let $D$ be a bounded convex domain in ${\Bbb C}^n$, $f:D\rightarrow D$ be a holomorphic endomorphism. Then if $Fix(f)$ is discrete and non-empty, it consists of one point only.
\eP
\begin{pf} Follows from the main theorem in [Vi1]: such a set has to be connected.
\end{pf}
Remark. Description of some properties of fixed point sets in convex domains can be found in [Ab].
\bT\label{bd} For any strictly pseudoconvex domain $D\subset {\Bbb C}^n$ with real analytic boundary, $n\ge 1$, the cardinality $\#(Fix(f))$ of the isolated fixed point set of an automorphism $f\in Aut(D)$ is finite. 
Moreover, there is a number $m=m(D)$ such that $\#(Fix(f))\leq m$. 
\eT
\begin{pf} If $D$ is biholomorphic to the ball or if $n=1$, then the statement is clear. Assume that $n\ge2$ and $D$ is not biholomorphic to the ball. By a theorem in [VEK], there is a neighborhood $U_1$ of $\overline D$ such that each automorphism of $D$ extends to be an injective holomorphic map on $U_1$. Consider a $g\in Aut(D)$. Choose domains $U_2, U_3$ with smooth boundaries so that $D\subset\subset U_3\subset\subset U_2\subset\subset U_1$. For every $h\in Aut(D)$ in some neighborhood of $g$, $h(\partial U_2)$ is so close to $g(\partial U_2)$ that $h(\partial U_2)\cap g(\overline U_3)=\emptyset$. Since $h(U_2)$ is a connected component of ${\Bbb C}^n\backslash h(\partial U_2)$ and since $h(U_2)\supset D$, we see that $h(U_2)\supset g(U_3)$ for every $h\in Aut(D)$ in some neighborhood of $g$. Since $Aut(D)$ is compact, there is a neighborhood $Q$ of $\overline D$ such that $Q\subset g(U_1)$ for each $g\in Aut(D)$. Let $U$ be the interior of the intersection of the sets $g(U_1)$, $g\in Aut(D)$. Then $U\supset Q$ and $g(U)=U$ for each $g\in Aut(D)$. I.e., each automorphism of $D$ is also an automorphism of $U$. There is
a finite cover of open sets $\{V_j: j=1, \dots, m\}$ of $\overline D$ such that each pair of points in a $V_j$ is connected by a unique distance-minimizing geodesic with respect to the Bergman metric of $U$. Let $f\in Aut(D)$. If $f$ fixes two points in a $V_j$, $f$ must fix each point on the unique
distance-minimizing geodesic connecting the two points. Consequently, each $V_j$ contains at most one isolated fixed point of $f$. Therefore, the number of isolated fixed points of $f$ is $\le m$.
\end{pf}

We now present counterexamples (e.g. the cardinality $\#(Fix(f))$ can be arbitrary, even infinity).
\bP For any $k\subset {\Bbb N}$, there exists a bounded domain $D\subset {\Bbb C}^n$, $n\ge 2$, and a holomorphic endomorphism $f:D\rightarrow D$, such that $\#(Fix(f))=k$.
\eP
\begin{pf} Without any loss of generality we can present an example for $n=2$. 

Let $S$ be the open Riemann surface in ${\Bbb C}^2$:
$$S=\{(x,y)\in C^2\mid y^2=(x-a_1)...(x-a_k)\},$$
where $a_1,...,a_k$ are $k$ distinct points in ${\Bbb C}$. The restriction $g$ of $(x,y)\mapsto (x,-y)$ to $S$ has exactly $k$ fixed points. Following [GR,
VIII, C8, p.257] there exists a holomorphic retraction $\rho :V\rightarrow S$
of an open neighborhood $V$ of $S$ onto $S$. Now the mapping $f:=g\circ \rho
:V\rightarrow V$ has exactly $k$ fixed points. Of course $V$ is not bounded,
but we can consider a bounded open set $W\subset V$ , $W\ni (a_s,0)$ for all 
$s=1,...,k$ and such that $g(W)=W$. This bounded domain will have the same
property.
\end{pf}

\bP There exists a hyperbolic manifold with a holomorphic automorphism whose fixed point set is discrete and consists of an infinite number of points.
\eP
\begin{pf} Consider the submanifold $X$ of $D^2$ defined by $y^2=B(x)$, where $D$ is the open unit disc and $B$ is a Blaschke product with an infinite number of zeroes, the restriction to $X$ of the map $(x,y)\rightarrow (x,-y)$ is an automorphism of $X$ and has an infinite number of isolated fixed points.
\end{pf}
\bP For any $n\ge 2$ and any $k\in {\Bbb N}$, there exists a polynomial automorphism $f$ of ${\Bbb C}^n$, such that $\#(Fix(f))=k$.
\eP
\begin{pf}. Let $a_1,...,a_k$ be k distinct complex numbers. Consider the map $H: {\Bbb C}^n \rightarrow {\Bbb C}^n$ given by
\newline $w_1 = z_1+z_2+(z_1-a_1)(z_1-a_2)...(z_1-a_k)$
\newline $w_2 = z_2+(z_1-a_1)(z_1-a_2)...(z_1-a_k)$
\newline $w_s=iz_s$ for all $s=3,...,n$
\newline One can easily check that this map is an automorphism [$(z_1,...,z_n)$ can be represented as polynomials of $(w_1,...,w_n)$], whose fixed point set is the set of the following $k$ points: $(a_1,0,...,0),(a_2,0,...,0),.....,(a_k,0,...,0)$.
\end{pf}
\bC Let $n\ge 2$; $p_1,p_2,...,p_k$ are $k$ distinct points in ${\Bbb C}^n$.
Then there exists a polynomial automorphism $g\in Aut({\Bbb C}^n)$ such that $
Fix(g)=\{p_1,p_2,...,p_k\}$.
\eC
\begin{pf} Let $p_j=(a_j,b_j),$ $a_j\in {\Bbb C},$ $b_j\in {{\Bbb C}^{n-1}}$. Without
any loss of generality we assume that the $a_j$'s are all distinct (in case
they are not, one can first use an invertible linear transformation of $
{\Bbb C}^n$ to achieve this). Consider the polynomial transformation $
F:w_1=z_1,w^{\prime }=z^{\prime }+f(z_1)$, where $f:{\Bbb C}\to {\Bbb C}
^{n-1}$ is the Lagrange interpolation polynomial map satisfying  $f(a_j)=b_j$
. Then $F(a_j,0)=p_j$,  $j=1,...,k$, and $F\in Aut({\Bbb C}^n)$. If $H\in
Aut({\Bbb C}^n)$ is the automorphism in the proof of the previous
proposition, then the automorphism $g=F\circ H\circ F^{-1}$ is such that $
Fix(g)=\{p_1,p_2,...,p_k\}$.
\end{pf}
\section{Single points as fixed point sets}

Here we consider some statements when the fixed point set of a holomorphic
automorphism is one point, specifically: which points can be a fixed point
set of an automorphism. First we show that if every point can be a fixed point set for an automorphism of a hyperbolic manifold then this manifold must be homogeneous. Second we provide an example when there are infinite number of points in the domain, each of which is a fixed point set for some holomorphic automorphism.

\bT If every point of a hyperbolic manifold $D$ is
a fixed point set for some holomorphic automorphism of $D$, then $D$ is a homogeneous manifold. 
\eT
For some concrete cases we have the following
\bC (A) If in the above theorem  $D\subset \subset {\Bbb C}^2$, then $D$ is biholomorphic to the unit ball ${\Bbb B}^2$ or the polydisc ${\Bbb U}^2$.

(B) If in the above theorem  $D\subset \subset {\Bbb C}^n$ has a smooth $C^2$
boundary, then $D$ is biholomorphic to the unit ball ${\Bbb B}^n$ .
\eC
To prove the Corollary we note that {\it (A)} there are only two kinds of bounded homogeneous domains in ${\Bbb C}^2$: the unit ball and the polydisc, and {\it (B)} in ${\Bbb C}^n$ there is only one bounded homogeneous domain with a smooth boundary: the unit ball (this is a consequence of a Wong-Rosay theorem (see [Ro], [Wo]). 
We will now prove the theorem.

\begin{pf} 
1. First we note that the theorem will follow from a local statement: let $
x\in D$, then there exists a neighborhood $U_x\,$of $x$ such that for any $
y\in U_x$ there is a $g\in Aut(D)$ such that $g(y)=x$.

Indeed, if this is true consider two arbitrary points $a,b\in D$, connect
them by a compact path $L$, cover $L$ by a finite number of $U_x$, $x\in L$,
and one can obtain an $f\in Aut(D)$, such that $f(a)=b$.

2. We now prove the local statement. Let $x\in D$.  By [FMV], for each point $x\in D$  there is an invariant Hermition metric in some neighborhood of the orbit $G(x)$, where $G=Aut(D)$. Consider a small enough
ball $b(x,\varepsilon )$ in that metric with center $x$ and radius $\varepsilon $,  
$\varepsilon >0$ will be determined by the construction later. Let $y\in
b(x,\varepsilon )$; consider the orbit $O(y)=\{z\in D:\exists g\in
Aut(D),g(y)=z\}$. Consider now a point $p\in O(y)$, such that $d(x,p)=d(x,O(y))$, where $d(\cdot ,\cdot )$ denotes the distance function induced by the local invariant metric.
Clearly, $p\in b(x,\varepsilon )$.  If $p=x$, there is nothing to prove;
otherwise consider a small ball $b_1$ of radius $<\frac
14d(x,p)$ that lies inside  $b(x,d(x,p))$, and such that  $\partial b_1\cap
\partial b(x,d(x,p) )=p$. 

This construction is possible if $\varepsilon $ is small enough, fixing such
an $\varepsilon =\varepsilon (x)$, we denote $b(x,\varepsilon )=U_x$. 

We observe that $O(y)\cap b(x,d(x,p))=\emptyset $. Let $q$ denote the center
of the ball $b_1$. By the assumption of the theorem there exists an $h\in
Aut(D)$ whose fixed point set is $q$. Now $h(p)\neq p$, and $h(p)\in
\partial b_1$, since $h(\partial b_1)=\partial b_1$. We now conclude that $
h(p)\in O(y)\cap b(x,d(x,p))$, which contradicts the previous observation
that this intersection is empty. Therefore $x=p\in O(y)$, and the theorem
has been proved.
\end{pf}

We now provide the following example.
\bP There exists a domain $D$ in ${\Bbb C}$ with infinite number of points each of which is the fixed point set for a holomorphic automorphism of $D$. 
\eP
\begin{pf} Consider $D={{\Bbb C}\backslash }\bigcup\limits_{n\in {\Bbb Z}}\Delta (n,1/3)$ where $\Delta (n,1/3)$ is a disk with center at $n\in {\Bbb Z}$ and radius $1/3$. Consider $f_k:z\mapsto (-z+(2k+1))$. Then for any $k\in {\Bbb Z}$ , $f_k\in Aut(D)$, and its fixed point set consists of one point $Fix(f)=\{k+1/2\}$.
\end{pf}
\section{Pairs of points as fixed point sets}

Here we examine the situation when the fixed point set consists of exactly
two points. Though such domains exist, no domain can have too
many pairs of distinct points as a fixed point set for an automorphism.

\bT Let $D\subset \subset {\Bbb C}^n$. The set $N\subset D^2\,$of all pairs,
each of which cannot be a fixed point set for a holomorphic automorphism of $D$, 
contains a full measure set in $D^2$.
\eT
It follows from the following
\bL\label{one} Let $D\subset \subset {\Bbb C}^n$, $a\in D$. Then there exists a complex
analytic set $Z\subset D$ ($\dim Z<n$ ), such that if $b\in D\backslash Z$
then the two points $\{a,b\}$ are such that for any automorphism $f
$ fixing these two points, the fixed point set of $f$ is at least one
(complex) dimensional.
\eL
First we need the following
\bL[H. Cartan]([Ca1, p.80])
{Let $D\subset \subset {\Bbb C}^n$, let $z\in D$, and let $I_z=I_z(D)$
be the isotropy subgroup at $z$ of the automorphism group of $D$. Then there
exists a holomorphic map $\phi :D\rightarrow {\Bbb C}^n$ such that $\phi
(z)=0,$  $\phi ^{\prime }(z)=id$, and for all $f\in I_z$
one has $\phi \circ f=f^{\prime }(z)\circ \phi $.}
\eL

\noindent  As in [Vi2, thm 2.3], for the proof of this Lemma, we define $\phi :D\to 
{\Bbb C}^n$ by 
\begin{equation*}
\phi (\zeta )=\int_{G_z}f^{\prime }(z)^{-1}(f(\zeta )-z)\,d\mu (f),
\end{equation*} 
where $d\mu $ is the Haar measure on $I_z$. Then $\phi (z)=0$, $\phi ^{\prime }(z)=id$ (and therefore $\phi$ is locally biholomorphic),
and $\phi \circ g=g^{\prime }(z)\circ \phi $ for each $g\in I_z$.

We are now ready to prove Lemma \ref{one}
\begin{pf}
Let $Z=\{z\in D|\varphi(z)=0\dot{\}}$. If $b\in D\backslash Z$, then suppose $f\in Aut(D)$ and $f$ fixes both points $a$ and $b$. We have $f^{\prime}(a)\cdot \varphi (b)=\varphi(f(b))=\varphi (b)$. Since by choice $\varphi (b)\neq 0$, and $\varphi $ is biholomorphic in the neighborhood $U$ of $a$, for a number $\lambda $, $|\lambda |>0$, small enough, there exists a point $c\in U\subset D$, $c\neq a$, $\varphi (c)=\lambda \varphi (b)$, and $f(c)\in U$. Now $\varphi(f(c))=f^{\prime}(a)\cdot \varphi (c)=f^{\prime}(a)\cdot \lambda \varphi (b)=\lambda \varphi (b)=\varphi (c)$.

Since $\varphi $ is biholomorphic in $U$ we have $f(c)=c.$
\end{pf}

\section{Application for holomorphic retractions}

Obviously any one point in a domain can be the fixed point set for an {\it endomorphism} of this domain. The theorem below gives an example of a domain that no two
distinct points can be the fixed point set for an endomorphism. Moreover:

\bT There is a domain $D\subset \subset {\Bbb C}^n$ such that for any two
distinct points $p\neq q\in D$, if a holomorphic endomorphism $f:D\rightarrow D$ fixes
these two points $(f(p)=p,f(q)=q)$ then $f=id$.
\eT
\begin{pf}
We will construct the example in ${\Bbb C}^2$; ${\Bbb C}^n$ with $n>2$ can be dealt with similarly.
	
We denote $B(z,r)$ - the euclidean ball with center at $z$ and of radius $r$, $b(z,r)\subset D\subset \subset {\Bbb C}^2$ - ball in the Kobayashi metric (in $D$) with center at $z\in D$ and of (Kobayashi) radius $r$. $B=B(0,1)={\Bbb B}^2$ the unit ball in ${\Bbb C}^2$.
	
1. Statement. {\it Let $a,b\in B$ be two distinct points, $L$ is the complex line through these points. Let $f\in H(B,B)$ fix these two points. Then $f$ fixes all the points of $L\cap B$.}
	
Proof of this statement follows from Example 1, Section 4 in [Vi3].
	
2. Statement. {\it If three distinct points $a,b,c\in B$ do not belong to the same complex line, then if $f\in H(B,B)$ fixes these points then $f=id $.}
\newline Proof of this statement follows from example 1, sec. 4 in [Vi3].

3. Statement. {\it $\forall a\in B_1=B(0,1/2)$, $a\neq 0$ there exists a $\it{unique}$ point $p\in \partial B_1$ that is closest to $a$ in the Kobayashi metric of the larger ball $B$: $k(a,p)=\underset{l\in \partial B_1}{\min k(a,l)}$, where $k(\cdot,\cdot)$ is the Kobayashi distance in $B$. Moreover, there exists a real number $s$ such that $p=s\cdot a$}.
\newline
To prove this let $\sigma $ be the Kobayashi distance from $a$ to $\partial B_1$. There exists an $r$ such that $b(0,\sigma )=B(0,r)$. Consider now $f\in Aut(B)$ such that $f(0)=a$. Then $f(B(0,r))=f(b(0,\sigma ))=b(a,\sigma )$, and from this construction we conclude that $\partial b(a,\sigma )$ and $\partial B_1\,$have only one common point $p$, moreover the vector $p$ is the intersection of the real line $\{s\cdot a\mid s\in {\Bbb R}\}$ with $\partial B_1$.
	
4. {\it Example of a domain $D\subset \subset {\Bbb C}^2$ and two points $\{a,b\}\in D$ such that any endomorphism of $D$ fixing those two points must be the identity.}
	
Consider $D=B\backslash \overline{B_1}$ and two distinct points $a,b\in D$ such that for the complex line $L\,$connecting them $L\cap B_1\neq \emptyset $ and $0\notin L$. Suppose $f\in H(D,D)$ fixes both points $a,b$. By Hartogs principle $f$ can be extended to $F\in H(B,B)$, and therefore $F$ fixes $L\cap B$. One can now pick a point $p\in L\cap B_1\,$ so that 
\newline (1) the boundary $\partial B_1$ has a unique point $c\in \partial B_1$ closest (in the Kobayashi metric of $B$) to $p$, and 
\newline (2) $c$ is not lying on $L$. 
\newline Since the Kobayashi metric cannot increase under holomorphic maps, $F(c)=c$. Now the three points $a,b,c\in B$ do not lie on the same complex line and by the previous statement $F$ (and therefore $f$) must be the identity.
	
5. We are now ready to prove the theorem by providing the main example: {\it a domain $D\subset \subset {\Bbb C}^2$ such that an endomorphism fixing any two given distinct points $\{a,b\}\in D$ is the identity.}
	
All we need to do is to take $B$ and remove a (countable) number of closed neighborhoods of portions of spheres, so that any complex line intersecting $B$ will intersect at least one of these removed spheres and then use the approach of the previous example.

For $k=1,2,\dots$, let
$$\varepsilon _k=1/2^{(4k)!},\;\;\;\delta _k=1/2^{(4k)!+1},$$
$$\Omega _{2k+1}=\overline{{B(0,1-\delta _k)\backslash B(0,1-\varepsilon _k)}}
\cap \{{Im}z_2\geq -1/2\},$$
and
$$\Omega _{2k+2}=\overline{B(0,1-4\delta _k)\backslash B(0,1-4\varepsilon _k)}\cap \{{Im}z_2\leq 1/2\}.$$
We will also need two more sets defined differently: 
$$\Omega _1=\overline{B(\alpha ,1-\delta _1)\backslash B(\alpha ,1-\varepsilon _1)}\cap \{{Im}z_2\geq -1/2\};$$
$$\Omega _2=\overline{B(\beta ,1-4\delta _1)\backslash B(\beta,1-4\varepsilon _1)}\cap \{{Im}z_2\leq 1/2\},$$
where $\alpha =(2^{-8!},0)$, $\beta =(0,2^{-8!})$. 
And finally $D=B\backslash (\bigcup\limits_{s=1}^\infty \Omega _s)$.
Now 	
$D$ is a connected open set. 
Let $a,b\in D$ and $f\in H(D,D)$ fixes $a,b$. Then $f$ can be extended to a holomorphic function $F:B\rightarrow B$. Let $L$ be the complex line connecting $a,b$, then (see Statement 1 above) $F|_L=id$.
	
$L\,$intersects an infinite number of $\Omega _s$. If $0\notin L$ there will always be at least one of two possibilities: either for some $s=2k+1$  there is a point $z\in L\cap \partial B(0,1-\delta _k)$, and ${Im}(z)>-1/2$, or for some $s=2k+2$ there is a point $z\in L\cap \partial B(0,1-4\delta _k)$, and ${Im}(z)<1/2$.
	
Similarly , if $0\in L$ there will always be at least one of two possibilities: either for $s=1$ there is a point $z\in L\cap \partial B(\alpha ,1-\delta _1)$, and ${Im}(z)>-1/2$, or for $s=2$ there is a point $z\in L\cap \partial B(\beta ,1-4\delta _1)$, and ${Im}(z)<1/2$. The above choice of $s$ is restricted in the following two cases: if $\alpha \in L$, we pick $s=2$, if $\beta \in L$, we pick $s=1$.
In any case we fix this point $z\in L\cap \partial \Omega _s$. If a point $p\in L$ is close enough to $z$, then the closest (in Kobayashi metric of the ball $B$) point to $p\,$ in the boundary $\partial \Omega _s$  is a unique point $c\in \partial \Omega _s$ that does not lie on $L$. Since $F(p)=p$, $F$ is a continuous non-increasing map in the Kobayashi metric, and $F(\overline{D})\subset \overline{D}$, we conclude that $F(c)=c$.
	
Now $F$ fixes three points in $B$ that do not lie on the same complex line, and therefore $F=id$, so $f=id$.
\end{pf}

A map $D\rightarrow D$ is a retraction, if $f\circ f=f$. A trivial
retraction is either a constant map, or the identity. 
\bC The domain in the above theorem has no non-trivial holomorphic retraction.
\eC
\section{Final remarks, unsolved problems}
\subsection{Some problems}

The main question that remains open is this.
1. Let $D$ be a bounded domain in ${\Bbb C}^n,$ $f\in Aut(D)$, and $Fix(f)$
is a discrete set. Can $\#(Fix(f))=\infty $?

If one considers the domain $D\subset $ ${\Bbb C}^n$ which is a direct
product of $n$ annuli, one can then find an $f\in Aut(D)$ with  $
\#(Fix(f))=2^n$. So, the next natural unsolved question is

2. Let $n\geq 2$, $D$ be a bounded domain in ${\Bbb C}^n\,$, with a piecewise smooth
boundary, $f\in Aut(D)$, and $Fix(f)$ is a set of isolated points. Can $
\#(Fix(f))\geq 2^n+1$?
(As noted earlier, for $n=1$ the answer is negative [PL]).

A more restricted version of the above question is a generalization of Theorem 1.2.
 
3. Is there a number $m$ such that for any strongly pseudoconvex domain $
D\subset \subset {\Bbb C}^n$, $\partial D\in C^\infty$, and $f\in Aut(D)$, if $Fix(f)$ is a set of
isolated points, then $\#(Fix(f))\leq m,$ where $m=m(n)$ (i.e. $m$ depends
on the dimension only)?

The next question, in case of a positive answer, would be a generalization
of Proposition 1.1.

4. Let $D$ be a bounded contractible domain in $
{\Bbb C}^n,$ $f\in Aut(D)$, and $Fix(f)$ is a non-empty set of isolated
points. Is $\#(Fix(f))=1$? 

\subsection{A (long) Remark. }

We now turn to a connection of this paper with the notion introduced and studied in papers [FK1, FK2, FM, KK, FMV, Vi2, Vi3]: {\it determining sets}.

Let $M$ be a complex manifold. Let $H(M,M)$ denote the set of holomorphic endomorphisms of $M$, and $Aut(M)$ the set of automorphisms of $M$.

\begin{Definition}
A set $K\subset M$ is called a determining subset of $M$ with respect to $
Aut(D)$ ($H(M,M)$ resp.) if, whenever $g$ is an automorphism (endomorphism
resp.) such that $Fix(g)\supseteq K$, then $Fix(g)=M$ (e.g. $g$ is the identity map of $M$).
\end{Definition}

One can now introduce a generalized notion of {\it quasi-determining set}
for a complex manifold $M$: 

\begin{Definition}
A set $K\subset M$ is called a quasi-determining subset of $M$ with respect
to $Aut(D)$ ($H(M,M)$ resp.) if, whenever $g$ is an automorphism
(endomorphism resp.) such that $Fix(g)\supseteq K$, then $K$ is a
proper subset of $Fix(g)$.
\end{Definition}
Another way to state this definition: A set $K\subset M$ is called a
quasi-determining subset of $M$ with respect to $Aut(D)$ ($H(M,M)$ resp.) if
it cannot be the fixed point set of any automorphism (endomorphism resp.) of $M
$.

There is an obvious reformulation of a number of results in our paper by
using this notion. For example, Proposition 1.1 means that any two points in
a convex domain form a quasi-determining set; Theorem 1.2 states that any $
m+1$ points in $D$ form a quasi-determining set, etc.

This definition obviously leads to a number of other open questions, which will 
be addressed in the future.

\end{document}